# On the Smoothness of the Solution to the Two-Dimensional Radiation Transfer Equation

**Dean Wang**

The Ohio State University
201 West 19th Avenue, Columbus, Ohio 43210

wang.12239@osu.edu

## ABSTRACT

In this paper, we deal with the differential properties of the scalar flux $\phi(x)$ defined over a two-dimensional bounded convex domain, as a solution to the integral radiation transfer equation. Estimates for the derivatives of $\phi(x)$ near the boundary of the domain are given based on Vainikko's regularity theorem. A numerical example is presented to demonstrate the implication of the solution smoothness on the convergence behavior of the Diamond Difference (DD) method.

## 1. INTRODUCTION

The one-group radiation transfer problem in a three-dimensional (3D) convex domain reads as follows: find a function $\phi: \overline{G} \times \Omega \to \mathbb{R}_+$ such as

$$\sum_{j=1}^{3} \Omega_j \frac{\partial \phi(x,\Omega)}{\partial x_j} + \sigma(x)\phi(x,\Omega) = \frac{\sigma_s}{4\pi} \int_{\Omega} s(x,\Omega,\Omega')\phi(x,\Omega')d\Omega' + f(x,\Omega) , \ x \in G , \quad (1)$$

$$\phi(x,\Omega) = \phi_{\partial G}(x,\Omega) , \quad x \in \partial G, \ \Omega \cdot \hat{n}(x) < 0 , \quad (2)$$

where $\Omega$ denotes the direction of radiation transfer, $\partial G$ is the boundary of the domain $G \subset \mathbb{R}^3$, $\sigma$ is the extinction coefficient (or macroscopic total cross section in neutron transport), $\sigma_s$ is the scattering coefficient (or macroscopic scattering cross section), $s$ is the phase function of scattering with $\int_{\Omega} s(x,\Omega,\Omega')d\Omega' = 4\pi$, $f$ is the external source function, and $\hat{n}$ is the unit normal vector of the domain surface. Note that $\sigma_s(x) \leq \sigma(x)$ under the subcritical condition, and $s(x,\Omega,\Omega') = s(x,\Omega',\Omega)$.

Assuming the isotropic scattering $s(x,\Omega,\Omega') = 1$, isotropic source $f(x,\Omega) = \frac{f(x)}{4\pi}$, and isotropic boundary condition $\phi_{\partial G}(x,\Omega) = \frac{\phi_{\partial G}(x)}{4\pi}$, we can obtain the so-called Peierls integral equation of radiation transfer for the scalar flux $\phi(x)$ as follows:

$$\phi(x) = \frac{1}{4\pi} \int_G \frac{\sigma_s(y) e^{-\tau(x,y)}}{|x-y|^2} \phi(y) dy + \frac{1}{4\pi} \int_G \frac{e^{-\tau(x,y)}}{|x-y|^2} f(y) dy$$

$$+ \frac{1}{4\pi} \int_{\partial G} \frac{e^{-\tau(x,y)}}{|x-y|^2} \left| \frac{x-y}{|x-y|} \cdot n(x) \right| \phi_{\partial G}(y) dS_y \ , \tag{3}$$

$$\tau(x,y) = \int_0^{|x-y|} \sigma(r - \xi\Omega) d\xi \ , \tag{4}$$

where $dS$ is the differential element of the domain surface, $\tau(x,y)$ is the optical path between $x$ and $y$. One can find detailed derivation in [1].

For simplicity, we assume $\sigma$ and $\sigma_s$ are constant over the domain. Then Eq. (3) can be simplified as

$$\phi(x) = \int_G K(x,y) \phi(y) dy + \frac{1}{\sigma_s} \int_G K(x,y) f(y) dy$$

$$+ \frac{1}{\sigma_s} \int_{\partial G} K(x,y) \left| \frac{x-y}{|x-y|} \cdot n(x) \right| \phi_{\partial G}(y) dS_y \ , \tag{5}$$

where the 3D radiation kernel is given as

$$K(x,y) = \frac{\sigma_s e^{-\sigma|x-y|}}{4\pi|x-y|^2} \ . \tag{6}$$

The boundary integral term in the above equation can produce singularities in the solution. We omit its discussion in this paper. In other words, here we only consider the problem with the vacuum boundary condition, i.e., $\phi_{\partial G} = 0$. Thus, Eq. (6) can be treated as the weakly integral equation of the second kind:

$$\phi(x) = \int_G K(x,y) \phi(y) dy + f(x) \ , \quad x \in G \ , \tag{7}$$

where $G \subset \mathbb{R}^n$, $n \geq 1$, is an open bounded domain and the kernel $K(x,y)$ is weakly singular, i.e., $|K(x,y)| \leq C|x-y|^{-\nu}$, $0 \leq \nu \leq n$. Weakly singular integral equations arise in many physical applications such as elliptic boundary problems and particle transport.

Since $K(x,y)$ has a singularity at $x = y$, the solution of a weakly integral equation is generally not a smooth function and its derivatives at the boundary would become unbounded from a certain order. There was extensive research on the smoothness (regularity) properties of the solutions to weakly integral equations [2,3], especially those early work in neutron transport theory done in the former Soviet Union [4,5]. It is believed that Vladimirov first proved that the scalar flux $\phi(x)$ possesses the property $|\phi(x+h) - \phi(x)| \sim h \log h$ for the one-group transport problem with isotropic scattering in a bounded domain [4]. Germogenova analyzed the local regularity of the angular flux $\phi(x, \Omega)$ in a neighborhood of the discontinuity interface and obtained an estimate of the first derivative, which has the singularity near the interface [5]. Pitkaranta derived a local singular resolution showing explicitly the behavior of $\phi(x)$ near the smooth portion of the boundary [6]. Vainikko introduced weighted spaces and obtained sharp estimates of pointwise derivatives near the smooth boundary for multidimensional weakly singular integral equations [1].

There exists some previous research work on the regularity of the integral radiation transfer solutions [7,8]. However, the 2D kernel treated in those studies bears no physical meaning since it does not model the scattering correctly. A physically relevant 2D kernel can be found in [9]. In this paper, we rederive the 2D kernel by directly integrating the 3D kernel with respect to the third dimension. We examine the differential properties of the new 2D kernel and provide estimates of pointwise derivatives of the scalar flux according to Vainikko's regularity theorem for the weakly integral equation of the second kind.

The remainder of the paper is organized as follows. In Sect. 2, we derive the 2D kernel for the integral radiation transfer equation. We examine the derivatives of the kernel and show that they satisfy the boundedness condition of Vainikko's regularity theorem in Sect. 3. Then the estimates of local regularity of the scalar flux near the boundary of the domain are given. Sect. 4 presents numerical results to demonstrate that the rate of convergence of a numerical method can be affected by the smoothness of the exact solution. Concluding remarks are given in Sect. 5.

## 2. TWO-DIMENSIONAL RADIATION TRANSFER EQUATION

In this section, we derive the 2D integral radiation transfer equation from its 3D form, Eqs. (5) and (6). In 3D, $dy = dy_1 dy_2 dy_3$ and $|x-y| = \sqrt{(x_1-y_1)^2 + (x_2-y_2)^2 + (x_3-y_3)^2}$. Let $\rho = \sqrt{(x_1-y_1)^2 + (x_2-y_2)^2}$, then $|x-y| = \sqrt{\rho^2 + (x_3-y_3)^2}$. In a 2D domain $G \subset \mathbb{R}^2$, the solution function $\phi(x)$ only depends on $x_1$ and $x_2$ in Cartesian coordinates. Therefore, we only need to find the 2D radiation kernel, which can be obtained by integrating out $y_3$ as follows:

$$K(x,y) = \int_{-\infty}^{\infty} \frac{\sigma_s e^{-\sigma|x-y|}}{4\pi|x-y|^2} \, dy_3$$

$$= \frac{\sigma_s}{4\pi} \int_{-\infty}^{\infty} \frac{e^{-\sigma\sqrt{\rho^2 + (x_3-y_3)^2}}}{\rho^2 + (x_3-y_3)^2} \, dy_3 \, . \tag{8}$$

To proceed, we introduce the variables $t = \sigma\sqrt{\rho^2 + (x_3-y_3)^2}$ and $z = y_3 - x_3$. Then we substitute $dy_3 = dz = \frac{t}{\sigma\sqrt{t^2 - \sigma^2\rho^2}} \, dt$ into the above equation to have

$$K(x,y) = \frac{\sigma_s}{4\pi} \int_{-\infty}^{\infty} \frac{e^{-t}}{\frac{t^2}{\sigma^2}} \, dz = \frac{\sigma_s}{2\pi} \int_0^{\infty} \frac{e^{-t}}{\frac{t^2}{\sigma^2}} \, dz$$

$$= \frac{\sigma_s}{2\pi} \int_{\sigma\rho}^{\infty} \frac{e^{-t}}{\frac{t^2}{\sigma^2}} \frac{t}{\sigma\sqrt{t^2 - \sigma^2\rho^2}} \, dt$$

$$= \frac{\sigma_s \sigma}{2\pi} \int_{\sigma\rho}^{\infty} \frac{e^{-t}}{t\sqrt{t^2 - \sigma^2\rho^2}} \, dt$$

$$= \frac{\sigma_s}{2\pi\rho} \int_0^1 \frac{e^{-\frac{\sigma\rho}{t}}}{\sqrt{1-t^2}} dt$$

$$= \frac{\sigma_s}{2\pi\rho} \int_1^\infty \frac{e^{-\sigma\rho t}}{\sqrt{1-t^2}} dt \ . \tag{9}$$

The above derivation of the 2D kernel is much simpler than the process given in [9], where the integral form is obtained by the projection of the particle flight path on the 2D plane. Note that the last two integrals are the first Bickley-Naylor function $Ki_1(r) = \int_0^1 \frac{e^{-\frac{r}{t}}}{\sqrt{1-t^2}} dt = \int_1^\infty \frac{e^{-rt}}{\sqrt{1-t^2}} dt$ [10,11].

Now we show that the 2D kernel $K(x, y)$ has a singularity at $\rho = 0$ (i.e., $x = y$) as follows. If $\rho$ is sufficiently small, then we can have $\sigma\rho < 1$ and

$$K(x, y) = \frac{\sigma_s}{2\pi\rho} \int_0^1 \frac{e^{-\frac{\sigma\rho}{t}}}{\sqrt{1-t^2}} dt > \frac{\sigma_s}{2\pi\rho} \int_0^1 \frac{e^{-\frac{1}{t}}}{\sqrt{1-t^2}} dt > \frac{0.3\sigma_s}{2\pi\rho} \ . \tag{10}$$

Thus, it can be seen that $K(x, y)$ becomes unbounded at $\rho = 0$.

By replacing the 3D kernel as defined by Eq. (6) with the above one, Eq. (5) becomes the 2D integral radiation transfer equation. Notice that the surface integral in the last term on the right-hand side of Eq. (5) should be replaced with a line integral in the 2D domain.

**Remark 2.1.** By simply assuming a 2D scattering, Johnson and Pitkaranta derived a 2D kernel, i.e., $K(x, y) = \frac{e^{-|x-y|}}{|x-y|}$ (where $\sigma = 1$), which is however physically incorrect since the scattering is essentially a 3D phenomenon [7]. Hennebach et al. also used the same 2D kernel for analyzing the radiation transfer solutions [8]. In addition, the integral equations in other geometries such as slab or sphere can be obtained by following the same approach, and they can be found in [12].

Applying Banach's fixed-point theorem, we can prove the existence and uniqueness of the solution in the 2D domain by showing that $\int_G K(x, y) dy$ is bounded below unity as follows.

$$\int_G K(x, y) dy = \int_G \left( \frac{\sigma_s \sigma}{2\pi} \int_{\sigma\rho}^\infty \frac{e^{-t}}{t\sqrt{t^2 - \sigma^2\rho^2}} dt \right) dy$$

$$= \frac{\sigma_s \sigma}{2\pi} \int_G dy \int_{\sigma\rho}^\infty \frac{e^{-t}}{t\sqrt{t^2 - \sigma^2\rho^2}} dt$$

$$= \frac{\sigma_s \sigma}{2\pi} \int_G \rho \, d\varphi \, d\rho \int_{\sigma\rho}^\infty \frac{e^{-t}}{t\sqrt{t^2 - \sigma^2\rho^2}} dt \ , \tag{11}$$

where $\varphi$ is the azimuthal angle. By extending the above bounded domain to the whole space, we have

$$\int_G K(x,y)dy < \frac{\sigma_s \sigma}{2\pi} \int_0^\infty 2\pi\rho d\rho \int_{\sigma\rho}^\infty \frac{e^{-t}}{t\sqrt{t^2 - \sigma^2\rho^2}} dt$$

$$= \sigma_s \sigma \int_0^\infty \rho d\rho \int_{\sigma\rho}^\infty \frac{e^{-t}}{t\sqrt{t^2 - \sigma^2\rho^2}} dt \ . \tag{12}$$

Denoting $\zeta = \sigma\rho$, Eq. (12) is simplified as

$$\int_G K(x,y)dy < \frac{\sigma_s}{\sigma} \int_0^\infty \zeta d\zeta \int_\zeta^\infty \frac{e^{-t}}{t\sqrt{t^2 - \zeta^2}} dt$$

$$= \frac{\sigma_s}{\sigma} \int_0^\infty \int_\zeta^\infty \frac{e^{-t}}{t} \frac{\zeta}{\sqrt{t^2 - \zeta^2}} dt d\zeta = \frac{\sigma_s}{\sigma} \int_0^\infty \frac{e^{-t}}{t} dt \int_0^t \frac{\zeta}{\sqrt{t^2 - \zeta^2}} d\zeta$$

$$= \frac{\sigma_s}{\sigma} \int_0^\infty e^{-t} dt$$

$$= \frac{\sigma_s}{\sigma} \leq 1 \ . \tag{13}$$

Notice we have changed the order of integration to solve the integral. It is apparent that for there exists a unique solution, the subcritical condition, $\sigma_s \leq \sigma$, must be satisfied.

**Remark 2.2**. A similar proof can be found for the 2D kernel in terms of the Bickley-Naylor function [13]. Existence of the unique solution to the general neutron transport equation in $L^1$ and $L^\infty$ has been long established in [14]. Corresponding results in $L^p$ for $1 \leq p \leq \infty$ can be found in [15-18].

## 3.   SMOOTHNESS OF THE SOLUTIONS

We first introduce Vainikko's regularity theorem [1], which provides a sharp characterization of singularities for the general weakly integral equation of the second kind. Then we analyze the differential properties of the 2D radiation kernel and show that the derivatives are properly bounded. Finally, Vainikko's theorem is used to give the estimates of pointwise derivatives of the radiation transport solution.

### 3.1.   Vainikko's Regularity Theorem

Before we state the theorem, we introduce the definition of weighted spaces $\mathbb{C}^{m,\nu}(G)$ [1].

**Weighted space $\mathbb{C}^{m,\nu}(G)$**. For a $\lambda \in \mathbb{R}$, introduce a weight function

$$w_\lambda = \begin{cases} 1 & , \ \lambda < 0 \\ (1 + |\log\varrho(x)|)^{-1} & , \ \lambda = 0 \ , \quad x \in G \\ \varrho(x)^\lambda & , \ \lambda > 0 \end{cases} \tag{14}$$

where $G \subset \mathbb{R}^n$ is an open bounded domain and $\varrho(x) = \inf_{y \in \partial G} |x - y|$ is the distance from $x$ to the boundary $\partial G$. Let $m \in \mathbb{N}$, $\nu \in \mathbb{R}$ and $\nu < n$. Define the space $\mathbb{C}^{m,\nu}(G)$ as the set of all $m$ times continuously differentiable functions $\phi \colon G \to \mathbb{R}$ such that

$$\|\phi\|_{m,\nu} = \sum_{\alpha \leq m} \sup_{x \in G} \left( w_{\alpha - (n-\nu)} |D^\alpha \phi(x)| \right) < \infty \,. \tag{15}$$

In other words, a $m$ times continuously differentiable function $\phi$ on $G$ belongs to $\mathbb{C}^{m,\nu}(G)$ if the growth of its derivatives near the boundary can be estimated as follows:

$$|D^\alpha \phi(x)| \leq c \begin{cases} 1 & , \ \alpha < n-\nu \\ 1 + |\log \varrho(x)| & , \ \alpha = n-\nu \,, \quad x \in G, \ \alpha \leq m \,, \\ \varrho(x)^{n-\nu-\alpha} & , \ \alpha > n-\nu \end{cases} \tag{16}$$

where $c$ is a constant. The space $\mathbb{C}^{m,\nu}(G)$, equipped with the norm $\|\cdot\|_{m,\nu}$, is a complete Banach space.

After defining the weighted space, we introduce the smoothness assumption about the kernel in the following form: the kernel $K(x,y)$ is $m$ times continuously differentiable on $(G \times G) \backslash \{x = y\}$ and there exists a real number $\nu \in (-\infty, n)$ such that the estimate

$$\left| D_x^\alpha D_{x+y}^\beta K(x,y) \right| \leq c \begin{cases} 1 & , \ \nu + \alpha < 0 \\ 1 + \left| \log |x-y| \right| & , \ \nu + \alpha = 0 \,, \quad x, y \in G \\ |x-y|^{-\nu-\alpha} & , \ \nu + \alpha > 0 \end{cases} \tag{17}$$

where

$$D_x^\alpha = \left( \frac{\partial}{\partial x_1} \right)^{\alpha_1} \cdots \left( \frac{\partial}{\partial x_n} \right)^{\alpha_n} \,, \tag{18}$$

$$D_{x+y}^\beta = \left( \frac{\partial}{\partial x_1} + \frac{\partial}{\partial y_1} \right)^{\beta_1} \cdots \left( \frac{\partial}{\partial x_n} + \frac{\partial}{\partial y_n} \right)^{\beta_n} \,, \tag{19}$$

holds for all multi-indices $\alpha = (\alpha_1, \cdots, \alpha_n) \in \mathbb{Z}_+^n$ and $\beta = (\beta_1, \cdots, \beta_n) \in \mathbb{Z}_+^n$ with $\alpha + \beta \leq m$. Here the following usual conventions are adopted: $\alpha = \alpha_1 + \cdots + \alpha_n$, and $|x| = \sqrt{x_1^2 + \cdots + x_n^2}$.

Now we present Vainikko's theorem in characterizing the regularity properties of a solution to the weakly integral equation of the second kind [1].

**Theorem 3.1**. Let $G \subset \mathbb{R}^n$ be an open bounded domain, $f \in \mathbb{C}^{m,\nu}(G)$ and let the kernel $K(x,y)$ satisfy the condition (17). If the integral equation (1) has a solution, $\phi \in L^\infty(G)$ then $\phi \in \mathbb{C}^{m,\nu}(G)$.

**Remark 3.1**. The solution does not improve its properties near the boundary $\partial G$, remaining only in $\mathbb{C}^{m,\nu}(G)$, even if $\partial G$ is of class $\mathbb{C}^\infty$, and $f \in \mathbb{C}^\infty(G)$. A proof can be found in [1]. More precisely, for any $n$ and $\nu$ ($\nu < n$) there are kernels $K(x,y)$ satisfying (17) and such that Eq. (1) is uniquely solvable and, for a suitable $f \in \mathbb{C}^\infty(G)$, the normal derivatives of order $k$ of the solution behave near $\partial G$ as $\log \varrho(x)$ if $k = n-\nu$, and as $\varrho(x)^{n-\nu-k}$ for $k > n-\nu$.

## 3.2. Smoothness of the Radiation Transfer Solution

To apply the results of **Theorem 3.1** to the 2D integral radiation transfer equation, we need to analyze the kernel $K(x, y)$ and show that it satisfies the condition (17), i.e., $\left| D_x^\alpha D_{x+y}^\beta K(x, y) \right| \leq c|x - y|^{-1-\alpha}$. We can simply set $|\beta| = 0$ without loss of generality for our problem.

$\boldsymbol{\alpha = 0}$:

$$K(x, y) = \frac{\sigma_s}{2\pi\rho} \int_0^1 \frac{e^{-\frac{\sigma\rho}{t}}}{\sqrt{1-t^2}} dt < \frac{\sigma_s}{2\pi\rho} \int_0^1 \frac{1}{\sqrt{1-t^2}} dt$$

$$= \frac{\sigma_s}{2\pi\rho} \frac{\pi}{2} = \frac{\sigma_s}{4\rho}$$

$$\leq c|x - y|^{-1} \,. \tag{20}$$

$\boldsymbol{\alpha > 1}$: Let $\zeta = \sigma\rho = \sigma|x - y| = \sigma\sqrt{(x_1-y_1)^2 + (x_2-y_2)^2}$, then $K(x, y) = \frac{\sigma_s}{2\pi\zeta} \int_0^1 \frac{e^{-\frac{\zeta}{t}}}{\sqrt{1-t^2}} dt$. By the chain rule, $D_x^\alpha$ can be written as

$$D_x^\alpha = \left(\frac{\partial}{\partial x_1}\right)^{\alpha_1} \left(\frac{\partial}{\partial x_2}\right)^{\alpha_2} = \left(\frac{\partial}{\partial\zeta}\frac{\partial\zeta}{\partial x_1}\right)^{\alpha_1} \left(\frac{\partial}{\partial\zeta}\frac{\partial\zeta}{\partial x_2}\right)^{\alpha_2} = \frac{\partial^\alpha}{\partial\zeta^\alpha} \left(\frac{\partial\zeta}{\partial x_1}\right)^{\alpha_1} \left(\frac{\partial\zeta}{\partial x_2}\right)^{\alpha_2} \,. \tag{21}$$

Then we have

$$|D_x^\alpha K(x, y)| = \left| \frac{\partial^\alpha}{\partial\zeta^\alpha} \left(\frac{\partial\zeta}{\partial x_1}\right)^{\alpha_1} \left(\frac{\partial\zeta}{\partial x_2}\right)^{\alpha_2} K(x, y) \right|$$

$$= \left| \left(\frac{\partial\zeta}{\partial x_1}\right)^{\alpha_1} \left(\frac{\partial\zeta}{\partial x_2}\right)^{\alpha_2} \right| \left| \frac{\partial^\alpha}{\partial\zeta^\alpha} K(x, y) \right| \,, \tag{22}$$

where

$$\left| \frac{\partial\zeta}{\partial x_1} \right| = \sigma \left| \frac{(x_1-y_1)}{\sqrt{(x_1-y_1)^2 + (x_2-y_2)^2}} \right| \leq \sigma \,, \tag{23a}$$

$$\left| \frac{\partial\zeta}{\partial x_2} \right| = \sigma \left| \frac{(x_2-y_2)}{\sqrt{(x_1-y_1)^2 + (x_2-y_2)^2}} \right| \leq \sigma \,. \tag{23b}$$

Substituting Eqs. (23a) and (23b) into (22), we obtain

$$|D_x^\alpha K(x, y)| \leq \sigma^\alpha \left| \frac{\partial^\alpha}{\partial\zeta^\alpha} K(x, y) \right| \,. \tag{24}$$

Apparently, we only need to find the upper bound of $\left| \frac{\partial^\alpha}{\partial\zeta^\alpha} K(x, y) \right| \leq c\zeta^{-(\alpha+1)}$, which is shown as follows.

Dean Wang

$$\left|\frac{\partial^\alpha}{\partial\zeta^\alpha}K(x,y)\right| = \left|\frac{\partial^\alpha}{\partial\zeta^\alpha}\left(\frac{\sigma_s}{2\pi\zeta}\int_0^1\frac{e^{-\frac{\zeta}{t}}}{\sqrt{1-t^2}}\,dt\right)\right| = \frac{\sigma_s}{2\pi}\left|\frac{\partial^\alpha}{\partial\zeta^\alpha}\left(\frac{1}{\zeta}\int_0^1\frac{e^{-\frac{\zeta}{t}}}{\sqrt{1-t^2}}\,dt\right)\right|$$

$$= \frac{\sigma_s}{2\pi}\left|\frac{(-1)^\alpha}{\zeta^{\alpha+1}}\int_0^1\frac{e^{-\frac{\zeta}{t}}}{\sqrt{1-t^2}}\,dt + \frac{1}{\zeta}\frac{\partial^\alpha}{\partial\zeta^\alpha}\left(\int_0^1\frac{e^{-\frac{\zeta}{t}}}{\sqrt{1-t^2}}\,dt\right)\right|$$

$$\leq \frac{\sigma_s}{2\pi}\frac{1}{\zeta^{\alpha+1}}\left|\int_0^1\frac{e^{-\frac{\zeta}{t}}}{\sqrt{1-t^2}}\,dt\right| + \frac{\sigma_s}{2\pi}\frac{1}{\zeta}\left|\frac{\partial^\alpha}{\partial\zeta^\alpha}\left(\int_0^1\frac{e^{-\frac{\zeta}{t}}}{\sqrt{1-t^2}}\,dt\right)\right|$$

$$< \frac{\sigma_s}{2\pi}\frac{1}{\zeta^{\alpha+1}}\left|\int_0^1\frac{1}{\sqrt{1-t^2}}\,dt\right| + \frac{\sigma_s}{2\pi}\frac{1}{\zeta}\left|\frac{\partial^\alpha}{\partial\zeta^\alpha}\left(\int_0^1\frac{e^{-\frac{\zeta}{t}}}{\sqrt{1-t^2}}\,dt\right)\right|$$

$$= \frac{\sigma_s}{4}\frac{1}{\zeta^{\alpha+1}} + \frac{\sigma_s}{2\pi}\frac{1}{\zeta}\left|\frac{\partial^\alpha}{\partial\zeta^\alpha}\left(\int_0^1\frac{e^{-\frac{\zeta}{t}}}{\sqrt{1-t^2}}\,dt\right)\right|. \tag{25}$$

Now we need to find the bound for $\left|\frac{\partial^\alpha}{\partial\zeta^\alpha}\left(\int_0^1\frac{e^{-\frac{\zeta}{t}}}{\sqrt{1-t^2}}\,dt\right)\right|$ in the above equation, which is analyzed as follows.

$$\left|\frac{\partial^\alpha}{\partial\zeta^\alpha}\left(\int_0^1\frac{e^{-\frac{\zeta}{t}}}{\sqrt{1-t^2}}\,dt\right)\right| = \left|\int_0^1\frac{\frac{\partial^\alpha}{\partial\zeta^\alpha}\left(e^{-\frac{\zeta}{t}}\right)}{\sqrt{1-t^2}}\,dt\right| = \left|\int_0^1\frac{e^{-\frac{\zeta}{t}}}{t^\alpha\sqrt{1-t^2}}\,dt\right| = \left|\int_0^\zeta\frac{e^{-\frac{\zeta}{t}}}{t^\alpha\sqrt{1-t^2}}\,dt + \int_\zeta^1\frac{e^{-\frac{\zeta}{t}}}{t^\alpha\sqrt{1-t^2}}\,dt\right|$$

$$< \left|\int_0^\zeta\frac{e^{-\frac{\zeta}{t}}}{t^\alpha\sqrt{1-t^2}}\,dt + \frac{1}{\zeta^\alpha}\int_\zeta^1\frac{e^{-\frac{\zeta}{t}}}{\sqrt{1-t^2}}\,dt\right| < \left|\int_0^\zeta\frac{e^{-\frac{\zeta}{t}}}{t^\alpha\sqrt{1-t^2}}\,dt\right| + \frac{1}{\zeta^\alpha}\left|\int_0^1\frac{1}{\sqrt{1-t^2}}\,dt\right|$$

$$= \left|\int_0^\zeta\frac{e^{-\frac{\zeta}{t}}}{t^\alpha\sqrt{1-t^2}}\,dt\right| + \frac{\pi}{2\zeta^\alpha}$$

$$= \frac{1}{\zeta^\alpha}\left|\int_1^\infty\frac{t^\alpha e^{-t}}{\sqrt{\left(\frac{t}{\zeta}\right)^2-1}}\,dt\right| + \frac{\pi}{2\zeta^\alpha}$$

$$< c\zeta^{-\alpha}. \tag{26}$$

Note that in the above equation, $0 < \zeta < 1$ and $\left|\int_1^\infty\frac{t^\alpha e^{-t}}{\sqrt{\left(\frac{t}{\zeta}\right)^2-1}}\,dt\right| < \infty$, i.e., it is bounded. Substituting Eq. (26) into (25) and noting that $\zeta = \sigma|x-y|$, we obtain

$$|D_x^\alpha K(x,y)| < c|x-y|^{-\alpha-1}. \tag{27}$$

Finally, we conclude that the 2D radiation kernel satisfies the condition (17). Therefore, by **Theorem 3.1**, the estimates of derivatives of the scalar flux $\phi(x)$ for radiation transfer are the same as those for the general weakly integral equation of the second kind:

$$|D^\alpha \phi(x)| \le c \begin{cases} 1 & , \ \alpha < 1 \\ 1 + |\log\varrho(x)| & , \ \alpha = 1 \\ \varrho(x)^{1-\alpha} & , \ \alpha > 1 \end{cases}, \quad x \in \mathrm{G}. \tag{28}$$

**Remark 3.2**. The first derivative of the solution $\phi(x)$ behaves as $\log\varrho(x)$ and becomes unbounded as approaching the boundary. The derivatives of order $k$ behave as $\varrho(x)^{1-k}$ for $k > 1$. These pointwise estimates cannot be improved by adding more strong smoothness on the data and domain boundary.

**Remark 3.3.** We point out that the lack of smoothness in the exact solution could adversely affect the convergence rate of spatial discretization schemes for solving the radiation transfer equation [19-21]. The spatial discretization error of a numerical method can be expressed as $\left|\phi_j - \phi_j^N\right| \le Ch_j^p \left\|\phi^{(k)}\right\|_\infty$, where $\phi_j$ is the exact solution at cell $j$, $\phi_j^N$ is its numerical result, $h_j$ is the mesh size, $p$ is the order of accuracy, and $\phi^{(k)}$ is the $k$-th derivative of the scalar flux. According to the above regularity results, we have $\left\|\phi^{(k)}\right\|_\infty \le ch_j^{1-k}$, then $\left|\phi_j - \phi_j^N\right| \le Ch_j^{p+1-k}$. However, for any numerical methods, we usually have $k \ge p$, which means that the spatial convergence rate should asymptotically reduce to the first order of accuracy or even worse as the mesh is refined.

## 4. NUMERICAL RESULTS

In this section, we demonstrate how the regularity of the exact solution will impact the numerical convergence rate by solving the S$_N$ neutron transport equation in its original integro-differential form, using the classic second-order Diamond Difference (DD) method. The model problem is a 1cm × 1cm square with the vacuum boundary condition. Thus, there will be no complication from the boundary condition. The S12 level-symmetric quadrature set is used for angular discretization.

We analyze the following four cases: **Case 1**: $\sigma = 1$, $\sigma_s = 0$; **Case 2**: $\sigma = 1$, $\sigma_s = 0.8$; **Case 3**: $\sigma = 10$, $\sigma_s = 0$, and **Case 4**: $\sigma = 10$, $\sigma_s = 0.9$. For all the cases, the external source $f = 1$, is infinitely differentiable, i.e., $f \in \mathbb{C}^\infty(\mathrm{G})$. Cases 1 and 3 are pure absorption problems, while Case 3 is optically thicker. It is worth noting that the solutions are only determined by the external source for these two cases. Cases 2 and 4 include the scattering effects, while Case 4 is optically thicker and more diffusive. Both the scattering and external source contribute to the solution. The flux $L^1$ errors as a function of mesh size and the rates of convergence are summarized in Table I. The error distributions on the mesh 160 × 160 are plotted in Fig. 1. The reference solutions for all the cases are obtained on a very fine mesh, 5120 × 5120.

**Table I.** Scalar Flux $L^1$ errors and convergence rates.

| Mesh | Case 1 | | Case 2 | | Case 3 | | Case 4 | |
|---|---|---|---|---|---|---|---|---|
| ($N \times N$) | Error[1] | Rate[2] | Error | Rate | Error | Rate | Error | Rate |
| $10 \times 10$ | 2.87E-03 | | 3.59E-03 | | 2.31E-03 | | 9.29E-03 | |
| $20 \times 20$ | 7.95E-04 | 1.85 | 1.01E-03 | 1.83 | 8.12E-04 | 1.51 | 2.56E-03 | 1.86 |
| $40 \times 40$ | 2.90E-04 | 1.45 | 3.73E-04 | 1.44 | 2.31E-04 | 1.82 | 5.89E-04 | 2.12 |
| $80 \times 80$ | 1.14E-04 | 1.35 | 1.44E-04 | 1.37 | 5.19E-05 | 2.15 | 1.37E-04 | 2.10 |
| $160 \times 160$ | 5.04E-05 | 1.17 | 6.32E-05 | 1.19 | 1.32E-05 | 1.97 | 3.53E-05 | 1.96 |
| $320 \times 320$ | 2.46E-05 | 1.03 | 3.06E-05 | 1.04 | 3.61E-06 | 1.87 | 9.39E-06 | 1.91 |
| $640 \times 640$ | 1.31E-05 | 0.91 | 1.63E-05 | 0.91 | 1.11E-06 | 1.71 | 2.70E-06 | 1.80 |
| $1280 \times 1280$ | 6.26E-06 | 1.07 | 7.76E-06 | 1.07 | 3.87E-07 | 1.51 | 8.51E-07 | 1.66 |

1.  The $L^1$ error on each mesh grid is defined by averaging the absolute differences between the numerical scalar flux and reference solution on each cell.
2.  The convergence rate $= \ln\left(\frac{\text{Error}(N \times N)}{\text{Error}(2N \times 2N)}\right)/\ln 2$.

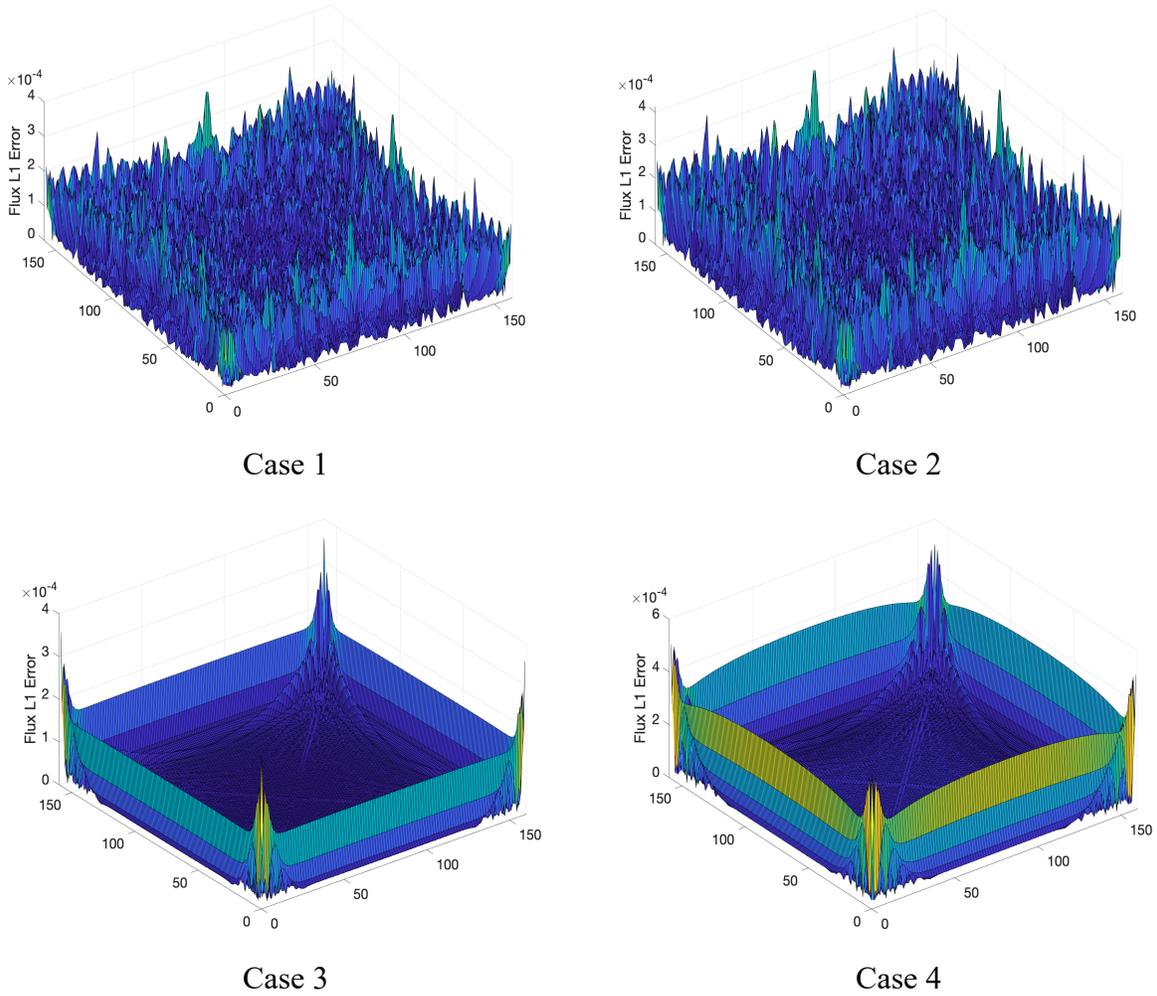

Case 1          Case 2

Case 3          Case 4

**Figure 1. Flux error distribution on the mesh $160 \times 160$.**

It is evident that the convergence rate decreases as the mesh is refined, and the errors are much larger at the boundary. The "noisier" distributions in Cases 1 and 2 are due to the ray effects of the discrete ordinates ($S_N$) method, which are more pronounced in the optically thin problem. The convergence behavior is similar between the two cases with and without the scattering, indicating that the source term plays a significant role in defining the irregularity of the solution. Cases 3 and 4 show the improved convergence rate as compared to Cases 1 and 2 because the exponential function $e^{-\sigma|x-y|}$ makes the kernel less singular for the larger total cross section $\sigma$. In addition, Case 4 has a slightly better rate of convergence than Case 3 on fine meshes (e.g., 1.80 vs. 1.71 on $640 \times 640$), because the transport problem in Case 4 becomes more like an elliptic diffusion problem [22], which usually has better regularity. It should be pointed out that in Case 3 the convergence rate is only 1.51 on the coarse mesh. It is because for the pure absorption case, the DD method becomes unstable when the mesh size is larger than $\frac{2\mu_j}{\sigma}$, where $\mu_j$ is the direction cosine of the radiation transfer direction. However, the stability condition can be somewhat relaxed for the scattering case.

**Remark 4.1**. The optimal error for the scalar flux of the DD method can be estimated by $|\phi_j - \phi_j^N| \le Ch_j^2 \|\phi''\|_\infty$ based on the discrete maximum principle [23]. For the angular flux, the truncation error estimates can be found in [24]. Although this optimal error estimate is obtained for the 1D slab geometry, one can expect the same to be true in two dimensions. As given by Eq. (28), the second derivative $\phi''$ will be bounded in the interior of the domain, while it would behave as $\phi'' \sim O(h_j^{-1})$ near the boundary. Therefore, it is expected that the convergence rate of the DD would decrease with refining the mesh, and asymptotically tend to $O(h_j)$. If the solution is sufficiently smooth (e.g., a manufactured smooth solution), the DD would maintain its second order of accuracy on any mesh size [25].

**Remark 4.2**. The scattering does not appear to play a role in defining the smoothness of the solution. For the problem without the external source, if there exists a nonsmooth or anisotropic incoming flux on the boundary, the scattering may not be able to regularize the solution either, since the irregularity caused by the incoming flux, which is defined by the surface integral term of Eq. (3), has nothing to do with the scattering and the solution flux $\phi$.

## 5. CONCLUSIONS

We have derived the two-dimensional integral radiation transfer equation and examined the differential properties of the integral kernel for fulfilling the boundedness conditions of Vainikko's theorem. We use the theorem to estimate the derivatives of the radiation transfer solution near the boundary of the domain. It is noted that the first derivative of the scalar flux $\phi(x)$ becomes unbounded when approaching the boundary. The derivatives of order $k$ behave as $\varrho(x)^{1-k}$ for $k > 1$, where $\varrho(x)$ is the distance to the boundary. A numerical example is presented to demonstrate that the irregularity of the exact solution will reduce the rate of convergence of numerical solutions. The convergence rate improves as the optically thickness of the problem increases. It is worth noting that the scattering does not help smoothen the solution. However, it does play a crucial role in transforming the transport problem into an elliptic diffusion problem in the asymptotic diffusion limit. We are currently extending the analysis to the boundary integral transport problem in considering nonzero incoming boundary conditions and corner effects. In addition, it would be interesting to study the convergence behavior of weak numerical solutions.